\documentclass{article}
\usepackage{amsmath,amsthm,amscd,amssymb,amsfonts}
\input{amssym.def}
\input{amssym.tex}

\newtheorem{thm}{Theorem}[section]
\newtheorem{prop}{Proposition}[section]
\newtheorem{lm}{Lemma}[section]

\newtheorem{remark}{Remark}[section]
\newtheorem{fact}{Fact}[section]
\newtheorem{cor}{Corollary}[section]

\title{An analog of the Iwasawa conjecture for a complete hyperbolic
threefold of finite volume}
\author{Ken-ichi SUGIYAMA
\footnote{Address : Ken-ichi SUGIYAMA,
 Department of Mathematics and Informatics,
Faculty of Science,
Chiba University,
1-33 Yayoi-cho Inage-ku,
Chiba 263-8522, Japan}
\footnote{e-mail address : sugiyama@math.s.chiba-u.ac.jp}}
\begin{document}
\maketitle
\begin{abstract}
For a unitary local system of rank one on a complete hyperbolic threefold of a
 finite volume with only one cusp, we will compare the
 order of its Alexander invariant at $t=1$ and one of the Ruelle L
 function at $s=0$. Our results may be considered as a
 solution of a geomeric analogue of the Iwasawa main conjecture in the algebraic number theory.
\end{abstract}
\section{Introduction}

In rescent days, it has been recognized there are many similarities
between the theory of a number field and one of a topological
threefold. In this note, we will show one more evidence, which is ``a
geometric analog of the Iwasawa conjecture''.\\

At first let us recall the original Iwasawa conjecture (\cite{Washington}). Let $p$ be an
odd prime and $K_n$ a cyclotomic field ${\mathbb Q}(\zeta_{p^n})$. The
Galois group ${\rm Gal}(K_n/{\mathbb Q})$ which is
isomorphic to ${\mathbb Z}/(p^{n-1})\times {\mathbb F}^{*}_p$ by the
cyclotomic character $\omega$ acts on
the $p$-primary part of the ideal class group $A_n$ of $K_n$. By the
action of ${\rm Gal}(K_1/{\mathbb Q})\simeq {\mathbb F}^{*}_p$, it has a
decomposition
\[
 A_{n}=\oplus_{i=0}^{p-2}A_{n}^{\omega^{i}},
\]
where we set
\[
 A_n^{\omega^i}=\{\alpha\in A_n \,|\, \gamma
      \alpha=\omega(\gamma)^i\alpha \,\,\mbox{for}\,\, \gamma \in {\rm
      Gal}(K_1/{\mathbb Q})\}.
\]
For each $i$ let us take the inverse limit with respect to the norm map:
\[
 X_i=\lim_{\leftarrow}A_n^{\omega^i}.
\]
If we set $K_{\infty}=\cup_{n}K_n$ and $\Lambda_{\infty}={\rm
Gal}(K_{\infty}/K_1)$, each $X_{i}$ becomes a ${\mathbb
Z}_{p}[[\Lambda_{\infty}]]$-module. Since there is an (noncanocal) isomorphism
${\mathbb Z}_{p}[[\Lambda_{\infty}]]\simeq {\mathbb Z}_{p}[[s]]$, each $X_i$ may be considered
as a ${\mathbb Z}_{p}[[s]]$-module. Iwasawa has shown that it is a
torsion ${\mathbb Z}_{p}[[s]]$-module and let ${\mathcal L}_{p}^{alg,i}$
be its generator, which will be referred as {\it the Iwasawa power
series.}\\

On the other hand, let 
\[
 {\mathbb Z}_{p}[[s]]\simeq {\mathbb Z}_{p}[[\Lambda_{\infty}]]\stackrel{\chi}\to{\mathbb Z}_{p}
\]
be the ring homomorphism induced by $\omega.$ For each $0 < i <p-1$, using the Kummer
congruence of the Bernoulli numbers, Kubota-Leopoldt and Iwasawa have
independently constructed an element of  ${\mathcal L}_p^{ana,i}$ which
satisfies
\[
 \chi^{r}({\mathcal L}_p^{ana,i})=(1-p^{r})\zeta(-r),
\]
for any positive integer $r$ which is congruent $i$ modulo $p-1$. Here
$\zeta$ is {\it the Riemann zeta function.} We will refer
${\mathcal L}_p^{ana,i}$ as {\it the $p$-adic zeta function}. The
Iwasawa main conjecture, which has been solved by Mazur and Wiles (\cite{MW})
says that ideals in ${\mathbb Z}_{p}[[s]]$ generated by ${\mathcal L}_{p}^{alg,i}$ and ${\mathcal
L}_p^{ana,i}$ are equal.\\

Now we will explain our geometric analog of the Iwasawa main
conjecture. \\

It is broadly recognized a geometric substitute for the
Iwasawa power series is the Alexander invariant. Let $X$ be a connected
finite CW-complex of dimension three and $\Gamma_{g}$ its fundamental
group. In what follows, we always assume that there is a surjective
homomorphism
\[
 \Gamma_{g} \stackrel{\epsilon}\to {\mathbb Z}.
\]
Let $X_{\infty}$ be the infinite cyclic covering of $X$ which
corresponds to ${\rm Ker}\,\epsilon$ by the geometric Galois theory and
$\rho$ a finite dimensional unitary representation of $\Gamma_g$. Then
$H_{\cdot}(X_{\infty},\,{\mathbb C})$ and $H_{\cdot}(X_{\infty},\,\rho)$
have an action of ${\rm
Gal}(X_\infty/X)\simeq {\mathbb Z}$, which make them $\Lambda$-modules.
Here we set $\Lambda={\mathbb C}[{\mathbb Z}]$ which is isomorphic to the Laurent polynomial ring ${\mathbb
C}[t,\,t^{-1}]$. Suppose that each of them is a torsion
$\Lambda$-module. Then due to the results of Milnor (\cite{MilnorI}), we
know $H^{i}(X_{\infty},\,\rho)$ is also a torsion $\Lambda$-module for
all $i$ and
vanishes for $i\geq 3$. Let $\tau^{*}$ be the action of $t$ on
$H^{i}(X_{\infty},\,\rho)$. Then {\it the Alexander invariant} is
defined to be the alternating product of the characteristic polynomials:
\[
 A_{\rho}^{*}(t)=\frac{\det[t-\tau^{*}\,|\, H^{0}(X_{\infty},\,\rho)]\cdot \det[t-\tau^{*}\,|\, H^{2}(X_{\infty},\,\rho)]}{\det[t-\tau^{*}\,|\, H^{1}(X_{\infty},\,\rho)]}.
\]

We will take {\it the Ruelle L-function} as a geometric substitute for
the $p$-adic zeta function. 
Let $\Gamma$ be a torsion free cofinite discrete subgroup of
$PSL_{2}({\mathbb C})$. It acts on the three dimensional Poincar\'{e}
upper half space
\[
 {\mathbb H}^3=\{(x,\,y,\,r)\,|\,x,\,y\in{\mathbb R},\,r>0\}
\]
endowed with a metric
\[
 ds^2=\frac{dx^2+dy^2+dr^2}{r^2},
\]
whose sectional curvature $\equiv -1$.
 Let $X$ be the quotient, which is a complete hyperbolic threefold of finite volume. We will
assume that it has only one cusp. 
Let $\rho$ be a unitary character of $\Gamma$. It defines a unitary
local system on $X$ of rank one, which will be denoted by the same symbol.
 By the one to one
correspondence between the set of loxiodromic conjugacy classes of
$\Gamma$ and one of closed geodesics of $X$,  {\it the Ruelle
L-function} is defined as
\[
 R_{\rho}(z)=\prod_{\gamma}\det[1-\rho(\gamma)e^{-zl(\gamma)}],
\]
where $\gamma$ runs through primitive closed geodesics.
Here $z$ is a complex number and $l(\gamma)$ is the length of
$\gamma$. 
It is known $R_{\rho}(z)$ is absolutely convergent if ${\rm Re}\,z$ is sufficiently large. We will show its square is
meromorphically continued on the whole plane. (If the
restriction $\rho|_{\Gamma_{\infty}}$ of $\rho$ to the fundamental group
$\Gamma_{\infty}$ of the cusp is nontrivial, $R_{\rho}(z)$ 
will be meromorphically
continued itself.) Let us define the order of $R_{\rho}(z)$
at $z=0$ to be
\[
 {\rm ord}_{z=0}R_{\rho}(z)=\frac{1}{2}{\rm ord}_{z=0}R_{\rho}(z)^{2}.
\]
We will compute it in terms of the dimension $h^{i}(\rho)$ of
$H^{i}(X,\,\rho)$.
\begin{thm}
Suppose $\rho|_{\Gamma_{\infty}}$ is trivial. Then we have
\[
 {\rm ord}_{z=0}R_{\rho}(z)=2(2h^{0}(\rho)-h^{1}(\rho)+1).
\]
On the contrary if $\rho|_{\Gamma_{\infty}}$ is nontrivial, 
\[
 {\rm ord}_{z=0}R_{\rho}(z)=-2h^{1}(\rho).
\]
\end{thm}
Note that in the latter case, $h^{0}(\rho)$ vanishes by the
assumption. We will find the ``error term'' $2$ in the first
identity is caused by the Hodge theory.\\

Suppose there is a surjective homomorphism from $\Gamma$ to ${\mathbb
Z}$ and let $X_{\infty}$ be the corresponding infinite covering of
$X$. Moreover suppose that all of the dimensions of
$H_{\cdot}(X_{\infty},\,{\mathbb C})$ and $H_{\cdot}(X_{\infty},\,\rho)$
are finite.
Let $g$ be a generator of the infinite cyclic group.
In \cite{Sugiyama} we have shown
\begin{thm}Suppose
that $H^{0}(X_{\infty},\,\rho)$ vanishes. Then
\[
 {\rm ord}_{t=1}A_{X}^{*}(\rho)\leq -h^{1}(\rho).
\]
Moreover if the action of $g$ on $H^{1}(X_{\infty},\,\rho)$ is
 semisimple, they are equal.
\end{thm}
These two theorems imply the following corollary, which may be
considered as a geometric analogue of the Iwasawa main conjecture.
\begin{cor} Suppose that $H^{0}(X_{\infty},\,\rho)$ vanishes.
\begin{enumerate}
\item If $\rho|_{\Gamma_{\infty}}$ is nontrivial, we have
\[
 {\rm ord}_{z=0}R_{\rho}(z)\geq 2{\rm ord}_{t=1}A_{X}^{*}(\rho).
\]
\item If $\rho|_{\Gamma_{\infty}}$ is trivial, we have
\[
 {\rm ord}_{z=0}R_{\rho}(z)\geq 2(1+{\rm ord}_{t=1}A_{X}^{*}(\rho)).
\]
\end{enumerate}
Moreover if the action of $g$ on $H^{1}(X_{\infty},\,\rho)$ is
 semisimple, they are equal.
\end{cor}
For a closed hyperbolic threefold, similar results have been proved in
\cite{Sugiyama}.\\

Although it seems curious there is a difference between two invariants,
such a phenomenon also occurs in the Iwasawa theory of an elliptic curve
(\cite{Bertolini-Darmon}). In their case the reason of the pathology is a
ramification and non-semisimplicity of a Galois representation
associated to an ellptic curve. This is quite similar to our case.\\

{\bf Acknowledgement.}\hspace{5mm}
It is a great pleasure to appreciate Professor Park for his kindness to
answer our many questions, as well as Professor Wakayama for sending his
manuscripts which were great hepl for us. It is clear without their help
our work will not be completed.

\section{A spectral decomposition and the Hodge theory}

Let $\Omega^j(\rho)$ be a vector bundle of $j$-forms on $X$
twisted by $\rho$ and the space of its square integrable sections will be
denoted by $L^2(X,\,\Omega^j(\rho))$. The positive Hodge Laplacian has the selfadjoint extension to $L^2(X,\,\Omega^j(\rho))$,
which will be denoted by $\Delta$. Note that the Hodge star operator
induces an isomorphism of Hilbert spaces:
\begin{equation}
 L^2(X,\,\Omega^j(\rho))\simeq L^2(X,\,\Omega^{3-j}(\rho)),\quad j=0,\,1,
\end{equation}
which commutes with $\Delta$. \\

Let $L^2(X,\,\Omega^j(\rho))_{d}$ be the closure of a subspace of
$L^2(X,\,\Omega^j(\rho))$ spanned by eigenvectors of the Laplacian
and $L^2(X,\,\Omega^j(\rho))_{c}$ its orthogonal complement. As we will
see in $\S4$, $L^2(X,\,\Omega^j(\rho))_{c}$  may
be $0$ according to the behavior of $\rho$ at the cusp. Also it is known
that $L^2(X,\,\Omega^j(\rho))_{c}$ is generated by the eigenpacket of
the Eisenstein series.(See $\S4$) The space
\[
 {\mathcal H}^{j}(\rho)=\{\varphi \in L^2(X,\,\Omega^j(\rho))\,|\,
 \Delta\varphi =0\},
\]
will be referred as {\it the space of harmonic forms}. It is
contained in $L^2(X,\,\Omega^j(\rho))_{d}$ and in particular its
dimension is finite. We will call its dimension  {\it $L^2$-Betti number}
 and write it by $\beta_j({\rho})_{(2)}$. Note that by (1) we have
\[
 \beta_j({\rho})_{(2)}=\beta_{3-j}({\rho})_{(2)},\quad j=0,\,1.
\]
We will explain the relation between the $L^2$-Betti numbers and the
topological one.\\

Let $h^{j}(\rho)$ (resp. $h^{j}_{c}(\rho)$) be the dimension of
$H^{j}(X,\,\rho)$ (resp. of the image of the compact supported
cohomology group $H^{j}_{c}(X,\,\rho)$ in $H^{j}(X,\,\rho)$). Then Zucker has shown (\cite{Zucker}, see also the
introduction of \cite{Mazzeo-Phillips}):
\begin{enumerate}
\item 
\[
 \beta_{0}(\rho)_{(2)}=h^{0}(\rho),
\]
\item 
\[
 \beta_{1}(\rho)_{(2)}=h^{1}_{c}(\rho).
\]
\end{enumerate}
For $A>0$ we set
\[
 {\mathbb H}^{3}_{A}=\{(x,\,y,\,r)\,|\, x,\,y\,\in{\mathbb R},\,r\leq A\}
\]
and let $X_A$ be its image by the natural projection:
\[
 {\mathbb H}^{3}\stackrel{\pi}\to X.
\]
If $A$ is sufficiently large, the complement
\[
 Y_A=X\setminus X_A
\]
is homeomorphic to a product of a two dimensional torus $T^2$ and an
open interval $(A,\,\infty)$. In particular $X_A$ is a deformation
retract of $X$ and we have a commutative diagram:
\[
 \begin{CD}
H^{1}(X_A,\partial X_A,\, \rho) @>>> H^{1}(X_A,\,\rho) \\
@VVV                                  @VVV \\
H^{1}_{c}(X,\,\rho) @>>> H^{1}(X,\,\rho).
\end{CD}
\]
Here the vertical arrows are isomorphisms.  The above morphism is
completed by the exact sequence:
\[
 \to H^{i-1}(\partial X_{A},\,\rho) \to H^{i}(X_A,\partial X_A,\, \rho)
 \to  H^{i}(X_A,\,\rho) \to H^{i}(\partial X_{A},\,\rho) \to.
\]
Note that we have
\[
 H^{i}(\partial X_{A},\,\rho)\simeq H^{i}(T^{2},\,\rho).
\]
\begin{lm}
Suppose $H^{0}(\partial X_{A},\,\rho)$ vanishes. Then we have
\[
 \beta_1(\rho)_{(2)}=h^{1}(\rho).
\]
\end{lm}
{\bf Proof.}
With the Poincar\'{e} duality, the assumption implies $H^{2}(\partial
X_{A},\,\rho)$ also vanishes. On the other hand the index theorem tells us 
\[
 \chi(\partial X_{A},\,\rho)=0.
\]
Thus we have 
\[
 H^{1}(\partial X_{A},\,\rho)=0
\]
and
\[
 H^{1}(X_A,\partial X_A,\, \rho)
 \to  H^{1}(X_A,\,\rho)
\]
is an isomorphism. Now the desired result will follow from Zucker's result.
\begin{flushright}
$\Box$
\end{flushright}
\begin{lm}
Suppose $H^{0}(\partial X_{A},\,\rho)$ does not vanish. Then we have
\[
 \beta_1(\rho)_{(2)}=h^{1}(\rho)-1.
\]
\end{lm}
{\bf Proof.} The assumption implies the restriction of $\rho$ to
$\partial X_A$ is trivial and 
\[
 H^{0}(\partial X_{A},\,\rho)=H^{2}(\partial X_{A},\,\rho)={\mathbb
 C},\quad H^{1}(\partial X_{A},\,\rho)={\mathbb C}^2.
\]
We first suppose that $H^{0}(X,\,\rho)$ vanishes. By the Poincar\'{e}
duality we know that $H^{3}(X_A,\partial X_A,\, \rho)$ also vanishes and have an exact sequence:
\begin{eqnarray*}
 0 & \to & {\mathbb C}  \to  H^{1}(X_A,\partial X_A,\, \rho) \to
  H^{1}(X_A,\, \rho) \\
& \to & {\mathbb C}^2 \to H^{2}(X_A,\partial X_A,\, \rho)\to H^{2}(X_A,\,
 \rho)\to {\mathbb C}\to 0.
\end{eqnarray*}
This exact sequence and the identity 
\[
 \dim H^j(X_A,\partial X_A,\, \rho)=h^{3-j}(\rho),\quad j=0,\,1.
\]
will imply
\[
 h^{1}(\rho)=h^{2}(\rho).
\]
In particular we have
\[
 0 \to {\mathbb C} \to H^{1}(X_A,\partial X_A,\, \rho) \to H^{1}(X_A,\,
 \rho) \to {\mathbb C} \to 0,
\]
which shows
\[
 \beta_1(\rho)_{(2)}=h^{1}(\rho)-1.
\]
Next suppose that $H^{0}(X,\,\rho)$ does not vanish. Then $\rho$ is the
trivial representation and the restriction
\[
  H^{0}(X_A,\,\rho) \to H^{0}(\partial X_{A},\,\rho)
\]
is an isomorphism. The Poincar\'{e} duality implies that the connecting homomorphism
\[
 H^{2}(\partial X_{A},\,\rho)\to H^{3}(X_A,\partial X_A,\,\rho)
\]
is also isomorphic and thus
\begin{eqnarray*}
0 &\to& H^{1}(X_A,\partial X_A,\,\rho) \to H^{1}(X_A,\,\rho)\to {\mathbb C}^2\\
&\to& H^{2}(X_A,\partial X_A,\,\rho) \to H^{2}(X_A,\,\rho)\to 0.
\end{eqnarray*}
is exact.
Now we will obtain the desired result by the same argument as before. 
\begin{flushright}
$\Box$
\end{flushright}
Here is an example. \\

Let $K$ be a hyperbolic knot in the three
dimensional sphere and $X$ its complement. By definition $X$ admits a
hyperbolic structure of a finite volume and it is known that the
dimension of $H^{1}(X,\,{\mathbb C})$ is one. Thus {\bf Lemma 2.2} shows
$\beta_1({\mathbb C})_{(2)}$ vanishes. 
\section{The derivative of the Laplace transform of the heat kernel}

For a function $f$ on $[0,\, \infty)$, we define its {\it derivative of the Laplace transform} as
\[
 L^{\prime}(f)(z)=2z\int_{0}^{\infty}e^{-tz^2}f(t)dt,
\]
if RHS is absolutely convergent. It is convenient to introduce more general
transformation:
\begin{equation}
 {\mathcal L}^{\prime}(f)(s,z)=2z\int_{0}^{\infty}t^{s-1}e^{-tz^2}f(t)dt.
\end{equation}
Suppose that the integral is absolutely convergent for $s$
and $z$ sufficiently large and that it is continued to a meromorphic function
on an open domain $U$ of ${\mathbb C}^2$ whose pole does not contain $\{(1,z)\,|\,z\in{\mathbb
C}\}$. Then we define
\[
 L^{\prime}(f)(z)={\mathcal L}^{\prime}(f)(s,z)|_{s=1},
\]
on $U\cap \{(1,z)\,|\,z\in{\mathbb C}\}.$
For example, let us take 
\[
 f(t)=t^{\nu},
\]
where $\nu$ is a half integer. Then the integral
\[
 2z\int_{0}^{\infty}t^{\nu+s-1}e^{-tz^2}dt
\]
is absolutely congergent for $z>0$ and $s>-\nu$ and is
computed as
\begin{eqnarray*}
{\mathcal L}^{\prime}(t^{\nu})(s,z)&=& 2z\int_{0}^{\infty}t^{\nu+s-1}e^{-tz^2}dt\\&=& 2z^{(1-2\nu)-2s}\Gamma(s+\nu).
\end{eqnarray*}
This is meromorphc function on 
$U=\{(s,\,z)\,|\,s,\,z\in{\mathbb C},\,-\pi< {\rm Im}\,z<\pi\}$ and we obtain
\begin{equation}
 L^{\prime}(t^{\nu})(z)=2z^{-(1+2\nu)}\Gamma(1+\nu).
\end{equation}
Note that the RHS is analytically continued as a rational function on the whole plane. 
\\

 We will use the following notation:
\[
 \delta_{0,\rho}(t)={\rm Trace}[e^{-t\Delta}\,|\,L^{2}(X,\,\Omega^{0}(\rho))_{d}],
\]
\[
  \delta_{1,\rho}(t)={\rm Trace}[e^{-t\Delta}\,|\,L^{2}(X,\,\Omega^{1}(\rho))_{d}]-\delta_{0,\rho}(t).
\]
Let $\Sigma_{j,\rho}$ be the set of eigenvalues of the Laplacian on
$L^{2}(X,\,\Omega^{j}(\rho))_{d}$. For sufficiently large $z\in{\mathbb
R}$, the derivatives of Laplace transforms of $e^{t}\delta_{0,\rho}(t)$
and $\delta_{1,\rho}(t)$ are well defined. We will study thier
properties.
\begin{lm} Let $L^{\prime}_{1}(z)$ be the derivative of the Laplace transform of
 $\delta_{1,\rho}(t)$. Then it is meromorphically continued on the whole
 plane and has only simple poles whose residues are integers.   Moreover it satisfies the following properties:
\begin{enumerate}
\item $L^{\prime}_1(z)$ satisfies a functional equation:
\[
 L^{\prime}_{1}(-z)=-L^{\prime}_{1}(z).
\]
\item The residue of $L^{\prime}_{1}(z)$ at $z=0$ is $2( \beta_1({\rho})_{(2)}- \beta_0({\rho})_{(2)})$.
\end{enumerate}

\end{lm}
{\bf Proof.}
For $\lambda \geq 0$ and $z>0$, we have
\begin{eqnarray*}
 L^{\prime}(e^{-t\lambda})(z)&=&2z\int^{\infty}_{0}e^{-t(z^2+\lambda)}dt\\
&=&\frac{2z}{z^2+\lambda}\\
&=& \frac{1}{z-i\lambda}+\frac{1}{z+i\lambda}.
\end{eqnarray*}
Thus the equation: 
\[
 \delta_{1,\rho}(t)=\sum_{\alpha\in \Sigma_{1,\rho}}e^{-t\alpha}-\sum_{\beta\in \Sigma_{0,\rho}}e^{-t\beta},
\]
shows
\begin{eqnarray*}
L^{\prime}_1(z)=L^{\prime}(\delta_{1,\rho})(z)&=& \frac{2(\beta_1({\rho})_{(2)}- \beta_0({\rho})_{(2)})}{z}\\
&+& \sum_{\alpha\in \Sigma_{1,\rho},\alpha >0}( \frac{1}{z-i\alpha}+\frac{1}{z+i\alpha})-\sum_{\alpha\in \Sigma_{1,\rho},\beta >0}( \frac{1}{z-i\beta}+\frac{1}{z+i\beta}),
\end{eqnarray*}
which implies desired results.
\begin{flushright}
$\Box$
\end{flushright}
The same computation will show
\begin{eqnarray*}
L^{\prime}(e^{t}\delta_{0,\rho})(z)&=& \beta_0(\rho)_{(2)}(\frac{1}{z-1}+\frac{1}{z+1})\\
&+& \sum_{\beta\in \Sigma_{0,\rho}, 0<\beta\leq 1}(\frac{1}{z-\sqrt{1-\beta}}+\frac{1}{z+\sqrt{1-\beta}})\\
&+&\sum_{\beta\in \Sigma_{0,\rho}, \beta>1}(\frac{1}{z-\sqrt{\beta-1}i}+\frac{1}{z+\sqrt{\beta-1}i}),
\end{eqnarray*}
which implies the following lemma.

\begin{lm} 
Let us put:
\[
 L^{\prime}_{0}(z)=L^{\prime}(e^{t}\delta_{0,\rho})(z-1).
\]
Then $L^{\prime}_{0}(z)$ is meromorphically continued on the whole plane and has only simple poles whose residue an integer.  Moreover it satisfies a functional equation:
\[
 L^{\prime}_0(1+z)=-L^{\prime}_0(1-z).
\]
and 
\[
 {\rm Res}_{z=0}L^{\prime}_{0}(z)={\rm Res}_{z=2}L^{\prime}_{0}(z)=\beta_0(\rho)_{(2)}.
\]
\end{lm}

\section{Selberg trace formula}

In this section, we will review {\it the Selberg trace
formula} following
\cite{Sarnak-Wakayama}. (See also \cite{EGM} and \cite{Park}.)\\
 
Let $A$ be a split Cartan subgoup of $G=PSL_{2}({\mathbb C})$. The Lie
algebras of $G$ and $A$ will be denoted by ${\frak G}$ and ${\frak A}$, respectively. The choice of $A$ determines a positive root $\alpha$ of ${\frak G}$ and let $H$ be an
element of ${\frak A}$ satisfying
\[
 \alpha(H)=1.
\]
If we exponentiate a linear isomorphism:
\[
 {\mathbb R} \stackrel{h}\to {\frak A},\quad h(t)=tH, 
\]
we know $A$ is isomorphic to the multiplicative group of positive real
numbers ${\mathbb R}_{+}$ and will identify them.

Let $K\simeq SO_{3}({\mathbb R})$ be the maximal compact
subgroup. According to the Iwasawa decompostion $G=KAN$ an element $g$
of $G$ can be written as
\[
 g=k(g)a(g)n(g).
\]
We put $r(g)=a(g)^{-1}$. \\

Let $M$ be the centralizer of $A$ in
$K$, which is isomorphic to $SO_{2}({\mathbb R})$. It determines a
paraboloic subgroup with a Langlands decomposition: 
\[
 P=MAN.
\]
Let $D_{K}$ or $D_{M}$ be the set of dominant integral forms on ${\frak
H}_{K}\otimes {\mathbb C}$, or on ${\frak m}\otimes {\mathbb
C}$ respectively. Here ${\frak H}_{K}$ (resp. ${\frak m}$) is a
Cartan subalgebra of $K$ (resp. the Lie algebra of $M$). The there is a
natural bijection between $D_{K}$ (resp. $D_{M}$) and the set of
nonnegative integers ${\mathbb Z}_{\geq 0}$ (resp. ${\mathbb Z}$). For
$\sigma \in D_{M}$ (resp. $\lambda \in D_{K}$), let ${\mathbb
C}(\sigma)$ (resp. $\tau_{\lambda}$) be the corresponding highest weight
representaion of $M$ (resp. $K$). Concretely $\tau_{\lambda}$ is the
standard representation of $K$ on the space of homogeneous polynomials
of three variables of degree $\lambda$. In particular the cotangent
bundle $\Omega^{1}_{{\mathbb H}^{3}}$ of ${\mathbb H}^{3}$ is a
homogeneous vector bundle associated to $\tau_{1}$:
\begin{equation}
 \Omega^{1}_{{\mathbb H}^{3}}=G\times_{K}\tau_{1}.
\end{equation} 
Let $\Gamma_{\infty}$ be the intersection of $\Gamma$ and $N$. 
Since we have assumed that $\Gamma$ has no elliptic element it coincides
with $\Gamma \cap P$. \\

Now we recall the principal and the Eisenstein series. Let ${\frak
A}_{\mathbb C}$ be the complexification of ${\frak A}$. For $(\sigma,\,s)
\in D_{M}\times {\frak A}^{*}_{\mathbb C}\simeq {\mathbb Z}\times
{\mathbb C}$, let ${\mathcal H}_{\sigma,s}$ be the Hilbert space of
Borel measurable functions on $G$ which satisfies
\[
 f(xman)=a^{-1-s}\sigma(m)f(x), \quad x\in G,\, m\in M,\,a\in A,
\] 
and
\[
 ||f||^{2}_{K}=\int_{K}|f(k)|^{2}dk < \infty.
\]
The integral is taken with respect to the Haar measure on $K$ of total
vulume one.

Now we will define an action of $G$ on ${\mathcal H}_{\sigma,s}$ to be
\[
 (\pi_{\sigma,s}(g)f)(x)=f(g^{-1}x),
\]
which will be referred as {\it the principal series representation.}
It is known $\pi_{\sigma,s}$ is unitary if and only if $s$ is pure
imaginary. The Cartan involution $w$ yields an isomorphism:
\[
 \pi_{\sigma,s}\simeq \pi_{-\sigma,-s}.
\]
In our case since any nonzero element of $D_{M}$ is unramified,
$\pi_{\sigma,s}$ is not isomorphic to $\pi_{\sigma,-s}$ for $\sigma\neq 0$.  
We set

\[
 \pi(\sigma,\,s)=
\left\{
\begin{array}{ccc}
\pi_{0,s}&\mbox{if}&\sigma=0\\
\pi_{\sigma,s}\oplus\pi_{-\sigma,s}&\mbox{if}&\sigma\neq 0
\end{array}
\right.
\]
and let ${\mathcal H}(\sigma,\,s)$ be its representation space.\\

 For a $K$-finite vector $\varphi_{\sigma}\in{\mathcal H}(\sigma,\,s)$, we associate {\it an Eisenstein series} $E(\varphi_{\sigma},\rho,s)$, which
is a function on $G$ defined as
\[
 E(\varphi_{\sigma},\rho,s)(x)=\sum_{\gamma\in \Gamma_{\infty}\backslash
 \Gamma}\rho(\gamma)^{-1}r(\gamma x)^{1+s}\varphi_{\sigma}(\gamma x).
\] 
It is known that an Eisenstein series safisfies the following
properties:
\begin{enumerate}
\item $E(\varphi_{\sigma},\rho,s)$ absolutely convergents to a $C^{\infty}$-function on
      $\{s\in{\mathbb C}\,|\, {\rm Re}\, s >1\}\times G$, which is
      holomorphic and real analytic in $s$ and $x$, respectively. Moreover it is meromorphically
      conitinued on the whole plane.
\item For $\gamma \in \Gamma$,
\[
 E(\varphi_{\sigma},\rho,s)(\gamma x)=\rho(\gamma)E(\varphi_{\sigma},\rho,s)(x).
\]
\item Let $Z$ be the center of the universal envelopping algebra of
      ${\frak G}$. Then $E(\varphi_{\sigma},\rho,s)$ is
      $Z$-finite.
\item Let us fix $x\in G$. Then the function on $K$ which is defined to be:
\[
 k\in K \to E(\varphi_{\sigma},\rho,s)(xk)
\]
is a $C^{\infty}$-function.
\item If $s$ is pure imaginary, $E(\varphi_{\sigma},\rho,s)$ is
      absolutely square
      integrable on $\Gamma\backslash G$. 
\end{enumerate}

Next we will consider the spectral decomposition of
$L^{2}(X,\,\Omega^{j}(\rho))$. For $m\in \Sigma_{j,\rho}$, let $e_{m}$
be the corresponding eigenvector. Let
$L_{\infty}$ be the torus $\Gamma_{\infty}\backslash {\mathbb C}$ and
$|L_{\infty}|$ its volume. For each $\sigma \in D_{M}$, let
$\varphi_{\sigma}\in {\mathcal H}_{\sigma,s}$ be a nonzero $K$-invariant vector. Note that
by definition
an element of ${\mathcal H}_{\sigma,s}$ is determined by its restriction
to $K$ and, by the Frobenius reciprocity law, $K$-invariant part of ${\mathcal H}_{\sigma,s}$ is
isomorphic to ${\mathbb C}(\sigma)$.\\

\begin{enumerate}
\item Suppose the restriction of $\rho$ to $\Gamma_{\infty}$ is trivial.
We will treat the spectral decomposition according to the case $j=0$ and $j=1$
      separately. 
\begin{enumerate}
\item The case of $j=0$.\\

Suppose $f\in L^{2}(X,\,\Omega^{0}(\rho))$ is contained in the domain of
$\Delta$. Then it has a spectral expansion:
\[
 f=\sum_{m\in\Sigma_{0,\rho}}(f,\,e_{m})e_{m}+\frac{1}{4\pi|L_{\infty}|}\int^{\infty}_{-\infty}(f,\,E(\varphi_{0},\rho,it))E(\varphi_{0},\rho,it)dt.
\]
\item The case of $j=1$.\\

Since the restriction of $\tau_{1}\otimes{\mathbb C}$ to $M$ has a 
      decomposition:
\[
 \tau_{1}\otimes{\mathbb C}|_{M}\simeq {\mathbb C}(-1)\oplus {\mathbb
      C}(0) \oplus {\mathbb C}(1), 
\] 
an element $f$ in the domain
      of the Laplacian has an expansion:
\[
  f=\sum_{m\in\Sigma_{1,\rho}}(f,\,e_{m})e_{m}+\frac{1}{4\pi|L_{\infty}|}\sum_{\sigma=-1}^{1}\int^{\infty}_{-\infty}(f,\,E(\varphi_{\sigma},\rho,it))E(\varphi_{\sigma},\rho,it)dt.
\]
\end{enumerate}
\item Suppose the restriction of $\rho$ to $\Gamma_{\infty}$ is
      nontrivial. Then we will know every element of
      $C^{\infty}(X,\,\Omega^{j}(\rho))$ is cuspidal. In fact let us
      choose $\gamma\in \Gamma_{\infty}$ so that $\rho(\gamma)\neq 1$.
      For $f\in C^{\infty}(X,\,\Omega^{j}(\rho))$ we have
\[
 \int_{L_{\infty}}f(x)dx=\int_{L_{\infty}}f(\gamma x)dx=\rho(\gamma)\int_{L_{\infty}}f(x)dx,
\]
which shows
\[
 \int_{L_{\infty}}f(x)dx=0.
\]
This implies that, in the spectral expansion, we do not have
      terms of Eisenstein series. Thus we have {\it the eigenfunction expansion}: 
\[
 f=\sum_{m\in\Lambda_{j,\rho}}(f,\,e_{m})e_{m},
\]
for $j=0$ and $1$.  
\end{enumerate}
Now we will explain the Selberg trace formula. We want to compute the
trace of heat kernel:
\[
 {\rm Tr}[e^{-t\Delta}\,|\, L^{2}(X,\,\Omega^{j}(\rho))].
\]
On the geometric side it is computed by {\it the orbital integrals}:
\[
 {\rm Tr}[e^{-t\Delta}\,|\, L^{2}(X,\,\Omega^{j}(\rho))]={\mathcal
 I}_{j}(t)+{\mathcal H}_{j}(t)+{\mathcal U}_{j}(t),
\] 
where ${\mathcal I}_{j}(t)$, ${\mathcal H}_{j}(t)$ and ${\mathcal
U}_{j}(t)$ are {\it the identity}, {\it the hyperbolic} and {\it the unipotent orbital integral},
respectively. Each term will be discussed in the following sections
separately. \\

On the other hand, accoding to the type of spectrum, we have an orthogonal
decomposition:
\[
 L^{2}(X,\,\Omega^{j}(\rho))=L^{2}(X,\,\Omega^{j}(\rho))_{d}\oplus L^{2}(X,\,\Omega^{j}(\rho))_{c}.
\]
It is known that ${\rm Tr}[e^{-t\Delta}\,|\,
L^{2}(X,\,\Omega^{j}(\rho))_{c}]$ is computed as
\[
 {\rm Tr}[e^{-t\Delta}\,|\, L^{2}(X,\,\Omega^{j}(\rho))_{c}]=-{\mathcal
 T}_{j}(t)-{\mathcal S}_{j}(t),
\]
where ${\mathcal T}_{j}(t)$ and ${\mathcal S}_{j}(t)$ are {\it the
threshold} and {\it the scattering term}, respectively. They are defined in
terms of the Fourier coefficients of the Eisenstein
series and we will also compute them in $\S8$. As we have seen, if the restriction of $\rho$ to
$\Gamma_{\infty}$ is nontrivial, they will not appear. Thus we have {\it
the Selberg trace formula}:
\begin{equation}
{\rm Tr}[e^{-t\Delta}\,|\, L^{2}(X,\,\Omega^{j}(\rho))_{d}]={\mathcal
 I}_{j}(t)+{\mathcal H}_{j}(t)+{\mathcal U}_{j}(t)+\delta_{\rho}({\mathcal
 T}_{j}(t)+{\mathcal S}_{j}(t)), 
\end{equation}

where 
\[
 \delta_{\rho}=
\left\{\begin{array}{ccc}
1 & if & \rho\,|\,{\Gamma_{\infty}}=1\\
0 & if & \rho\,|\,{\Gamma_{\infty}}\neq 1
\end{array}
\right.
\]

\section{Ruelle L-function and hyperbolic terms}

Let $\Gamma_h$ be the set of
conjugacy classes of loxidromic elements of $\Gamma$. Since there is a
natural bijection between closed geodesics of $X$ and $\Gamma_h$, we may
identify them. Thus an element $\gamma$ of $\Gamma_h$ is written as 
\[
 \gamma=\gamma_0^{\mu(\gamma)},
\]
where $\gamma_0$ is a primitive closed geodesic and $\mu(\gamma)$ is a
positive integer, which will be referred as {\it the multiplicity}. 
The length of $\gamma\in \Gamma_h$ will be denoted by $l(\gamma)$.
Let $\Gamma_{h,prim}$ be the set of primitive closed geodesics.\\

Using the Langlands decomposition, $\gamma\in \Gamma_h$ may be written as
\[
 g\gamma g^{-1}=m(\gamma)\cdot a(\gamma)\in MA
\]
for a certain $g\in G$. Here $m(\gamma)$ is nothing but the
holonomy of a pararell transformation along $\gamma$. Note that elements of $GL_2({\mathbb R})$:
\[
 A^{u}(\gamma)=e^{l(\gamma)}m(\gamma), \quad A^{s}(\gamma)=e^{-l(\gamma)}m(\gamma)
\]
describe an unstable or a stable action of the linear Poincar\'{e}
map, respectively.\\

 For $\gamma\in \Gamma_h$ we set
\[
 \Delta(\gamma)=\det [I_2-A^{s}(\gamma)]
\]
and
\[
 a_{0}(\gamma)=\frac{\rho(\gamma)\cdot l(\gamma_0)}{\Delta(\gamma)},\quad
 a_{1}(\gamma)=\frac{\rho(\gamma)\cdot {\rm Tr}\,[m(\gamma)]\cdot l(\gamma_0)}{\Delta(\gamma)}.
\]
Now {\bf Theorem 2} of \cite{Fried} shows the hyperbolic terms are
given by
\[
 {\mathcal H}_0(t)=H_{0}(t),\quad {\mathcal H}_1(t)=H_{0}(t)+H_{1}(t),
\]
where 
\[
 H_{0}(t)=\Sigma_{\gamma\in \Gamma_h}\frac{a_0(\gamma)}{\sqrt{4\pi
 t}}\exp [-(\frac{l(\gamma)^2}{4t}+t+l(\gamma))],
\]
and
\[
 H_{1}(t)=\Sigma_{\gamma\in \Gamma_h}\frac{a_1(\gamma)}{\sqrt{4\pi
 t}}\exp [-(\frac{l(\gamma)^2}{4t}+l(\gamma))].
\]
We will explain a relation between these hyperbolic terms and the
Ruelle L-function. \\

 For $j=0,\,1$ we set
\[
 S_j(z)=\exp [-\sum_{\gamma\in \Gamma_{h}}\frac{a_j(\gamma)}{l(\gamma)}e^{-zl(\gamma)}],
\] 
and let $Y_j(z)$ be its logarithmic derivative:
\[
 Y_j(z)=\Sigma_{\gamma\in \Gamma_{h}}a_j(\gamma)e^{-zl(\gamma)}.
\]
Then Fried has shown (\cite{Fried} p.532, the formula (RS)):
\begin{equation}
 R_{\rho}(z)=\frac{S_0(z)\cdot S_0(z+2)}{S_1(z+1)}.
\end{equation}

\begin{remark}
If we use terminologies of Fried(\cite{Fried} p.529), for a hyperbolic
 threefold we have
\[
 \sigma_{0}=\sigma_{2},
\]
which implies
\[
 S_0=S_2.
\]
\end{remark}
Let $L^{\prime}_{0,hyp}(z)$ and $L^{\prime}_{1,hyp}(z)$ denote $L^{\prime}(e^{t}H_0)(z-1)$ and
$L^{\prime}(H_1)(z)$, respectively. 
\begin{prop}
\begin{enumerate}
\item
\[
 L^{\prime}_{0,hyp}(z)=Y_0(z).
\]
\item
\[
 L^{\prime}_{1,hyp}(z)=Y_1(z+1).
\]
\end{enumerate}
\end{prop}
{\bf Proof.} Since each statement will proved by the same way,  we will only
prove the first. By analytic continuation we may assume that $z>1$. The equation: 
\[
 e^{t}H_{0}(t)=\Sigma_{\gamma\in \Gamma_h}\frac{a_0(\gamma)}{\sqrt{4\pi
 t}}\exp [-(\frac{l(\gamma)^2}{4t}+l(\gamma))],
\]
and
\[
 2z\int^{\infty}_{0}\frac{1}{\sqrt{4\pi t}}\exp[-tz^2-\frac{l(\gamma)^2}{4t}]dt=e^{-zl(\gamma)}
\]
implies
\begin{eqnarray*}
L^{\prime}(e^tH_0)(z)&=& \Sigma_{\gamma\in \Gamma_h} a_0(\gamma)e^{-l(\gamma)}\cdot 2z\int^{\infty}_{0} \frac{1}{\sqrt{4\pi
 t}}\exp [-(tz^2+\frac{l(\gamma)^2}{4t})]dt\\
&=& \Sigma_{\gamma\in \Gamma_h} a_0(\gamma)e^{-(z+1)l(\gamma)}\\
&=& Y_0(z+1).
\end{eqnarray*}

\begin{flushright}
$\Box$
\end{flushright}
Combining {\bf Proposition 5.1} with (6), we have obtained
\begin{equation}
 \frac{d}{dz}\log R_{\rho}(z)=L^{\prime}_{0,hyp}(z)-L^{\prime}_{1,hyp}(z)+L^{\prime}_{0,hyp}(z+2).
\end{equation}

\section{Identity terms}

The formula of the Planchrel measures (\cite{Knapp}) and {\bf Theorem
2} of \cite{Fried} implies
\[
 {\mathcal I}_0(t)=I_{0}(t),\quad {\mathcal I}_1(t)=I_{0}(t)+I_{1}(t),
\]
where 
\[
 I_{0}(t)=vol(X)\int^{\infty}_{-\infty}e^{-t(x^2+1)}x^2dx,
\]
and
\[
 I_{1}(t)=2vol(X)\int^{\infty}_{-\infty}e^{-tx^2}(x^2+1)dx.
\]
We will compute the derivative of the Laplace transforms of $e^{t}I_{0}$ and
$I_1$.\\

Taking a derivative of the identity:
\[
 \int^{\infty}_{-\infty}e^{-tx^2}dx=\frac{\sqrt{\pi}}{2}t^{-\frac{1}{2}}
\]
with respect to $t$, we obtain
\[
 \int^{\infty}_{-\infty}x^{2}e^{-tx^2}dx=\frac{\sqrt{\pi}}{4}t^{-\frac{3}{2}}.
\]
In particular we have
\[
 \int^{\infty}_{-\infty}e^{-tx^2}(1+x^2)dx=\frac{\sqrt{\pi}}{2}(t^{-\frac{1}{2}}+\frac{t^{-\frac{3}{2}}}{2}).
\]
The identity
\[
 \sqrt{\pi}=\Gamma(\frac{1}{2})=-\frac{1}{2}\Gamma(-\frac{1}{2}),
\]
and the equation (3) will show
\[
 L^{\prime}(I_1)(z)=-2\pi vol(X)(z^2-1).
\]

By the same computation we will see
\[
 L^{\prime}(e^tI_0)(z)=-\pi vol(X)z^2.
\]
Thus we have obtained the following proposition.
\begin{prop}
\begin{enumerate}
\item
\[
 L^{\prime}(e^tI_0)(z)=-\pi vol(X)z^2.
\]
\item
\[
 L^{\prime}(I_1)(z)=-2\pi vol(X)(z^2-1).
\]
\end{enumerate}
\end{prop}

\section{Unipotent terms}

We will recall the Osborne and Warner's
formula.(\cite{Osborne-Warner})\\

Since the nilpotent radical $N$ is diffeomorphic to its Lie algebra ${\frak N}$ by the exponential map:
\[
 {\frak n}\stackrel{exp}\to N,
\]
the Killing form induces a norm $||\cdot||$ on $N$. We define 
{\it the Epstein L-function} of $\rho$ to be
\[
 L(\rho,\,s)=\sum_{0\neq\gamma\in \Gamma_{\infty}}\rho(\gamma)||\gamma||^{-2(1+s)}.
\]
It is absolutely convergent for ${\rm Re}\, s>0$ and is meromorphically
continued on the whole plane. Let $\rho|_{\Gamma_{\infty}}$ be the
restriction of $\rho$ to $\Gamma_{\infty}$. Then it is known
\begin{enumerate}
\item If $\rho|_{\Gamma_{\infty}}$ is not trivial,
      it is an entire function.
\item If  $\rho|_{\Gamma_{\infty}}$ is trivial, it
      has a simple pole only at $s=0$.
\end{enumerate}
Let $R_{\rho}$ be the its residue at $s=0$ and we put
\[
 C_{\rho}=\lim_{s\to 0}\{L(\rho,\,s)-\frac{R_{\rho}}{s}\}.
\]
Following Osborne and Warner (\cite{Osborne-Warner}), for a function on
$G$ which belongs to a Schwartz space ${\mathcal C}^p(G) \quad(0<p<1)$, we will consider the functions:
\[
 T(f,s)=\frac{1}{2\pi}\int_{N}dn ||n||^{-2s}\int_{K}f(knk^{-1})dk,
\]
and
\[
 I(f,s)=2|L_{\infty}|L(\rho,\,s)T(f,s).
\]
It is known that $T(f,s)$ is regular at $s=0$. In p.297 of
\cite{Osborne-Warner} they have shown the unipotent orbital integral of $f$
is given by 
\[
 {\mathcal U}(f)=\lim_{s\to 0}\frac{d}{ds}[s\cdot I(f,s)].
\]
This implies the following corollary.
\begin{cor}
\begin{enumerate}
\item If $\rho|_{\Gamma_{\infty}}$ is nontrivial, we have
\[
 {\mathcal U}(f)=2|L_{\infty}|C_{\rho}T(f,0).
\]
\item If $\rho|_{\Gamma_{\infty}}$ is trivial, we have
\[
 {\mathcal U}(f)=2|L_{\infty}|\{R_{\rho}T^{\prime}(f,0)+C_{\rho}T(f,0)\}.
\]

\end{enumerate}
\end{cor}
We will apply the corollary to our heat kernel. \\

Let $o$ be a point of
${\mathbb H}^3$ defined by
\[
 o=(0,\,0,\,1),
\]
and for $x\in G$ we set
\[
 [x]=x\cdot o.
\]
Let $\tilde{K}_j(\cdot,\cdot,t)$ be the integral kernel of the heat
operator on $\Omega^{j}_{{\mathbb H}^3}(\rho)$. Since
$\Omega^{j}_{{\mathbb H}^3}$ is a homogeneous vector bundle:
\[
 \Omega^j_{{\mathbb H}^3}=G\times_{K}\tau_j,
\]
there is a $K$-biinvarinat function
\[
 K_{j,t}\in C^{\infty}(G, End(\tau_j))
\]
such that
\[
 K_{j,t}(x^{-1}\cdot y)=\tilde{K}_j([x],[y],t), \quad x,y\in G. 
\]
Its trace $k_{j,t}$ is contained in ${\mathcal C}^{p}(G)$ for a certain
$0<p<1$. We will define its {\it nonabelain Fourier transform} $\hat{k}_{j,t}$,
 which is a function on the set $D_{M}\times i{\frak A}^{*}\simeq
{\mathbb Z}\times i{\mathbb R}$, to be:
\[
 \hat{k}_{j,t}(\sigma, i\lambda)=Tr[\pi_{\sigma, i\lambda}(k_{j,t})].
\]
Then Fried has shown the following result(\cite{Fried}{\bf Lemma 1}).
\begin{fact}
\[
 \hat{k}_{j,t}(\sigma, i\lambda)=
\left\{
\begin{array}{ccc}
e^{-t\lambda^2} & if & \sigma=\pm 1, \,j=1\\
e^{-t(1+\lambda^2)} & if & \sigma=0,  \,j=0, \,1\\
0 & otherwise &
\end{array}
\right.
\]
\end{fact}
 This fact and the formula of Park (\cite{Park},p.12) will imply
\begin{equation}
 T(k_{0,t},0)=\frac{e^{-t}}{4\pi^2}\int_{-\infty}^{\infty}e^{-t\lambda^2}d\lambda
\end{equation}
and
\begin{equation}
 T(k_{1,t},0)=\frac{1}{2\pi^2}\int_{-\infty}^{\infty}e^{-t\lambda^2}d\lambda+\frac{e^{-t}}{4\pi^2}\int_{-\infty}^{\infty}e^{-t\lambda^2}d\lambda.
\end{equation}
We put
\[
 U_0(t)={\mathcal U}_0(t),\quad U_1(t)={\mathcal U}_1(t)-{\mathcal U}_0(t).
\]
The following proposition will follow from {\bf Corollary 7.1}.
\begin{prop}
Suppose $\rho|_{\Gamma_{\infty}}$ is nontrivial. Then we have
\[
 U_0(t)=\frac{|L_{\infty}|C_{\rho}}{2\pi^2}\cdot e^{-t}\int_{-\infty}^{\infty}e^{-t\lambda^2}d\lambda,
\]
and
\[
 U_{1}(t)=\frac{|L_{\infty}|C_{\rho}}{\pi^2}\cdot\int_{-\infty}^{\infty}e^{-t\lambda^2}d\lambda.
\]
\end{prop}
We will compute the derivative of the Laplace transform of $e^tU_0$ and $U_1$.\\

For a nonzero real number $r$, let ${\rm sgn}(r)$ be its sign.
The following lemma will be proved by the Cauchy's integral formula.
\begin{lm}
Let $z$ be a positive number and $\alpha$ a complex number.
\begin{enumerate}
\item Suppose ${\rm Im}\, \alpha\neq 0$. Then we have
\[
 \int_{-\infty}^{\infty}\frac{d\lambda}{(\lambda^2+z^2)(\lambda-\alpha)}=\frac{\pi i}{z(z-{\rm sgn}({\rm Im}\,\alpha)i\alpha)}.
\]
\item 
\[
 \int_{-\infty}^{\infty}\frac{d\lambda}{\lambda^2+z^2}=\frac{\pi}{z}.
\]
\end{enumerate}
\end{lm}
For $z>0$, using {\bf Lemma 7.1}, we have
\begin{eqnarray*}
L^{\prime}(e^{t}U_{0})(z)&=&\frac{|L_{\infty}|C_{\rho}}{2\pi^2}\cdot 2z\int^{\infty}_{0}dt e^{-tz^2}\int_{-\infty}^{\infty}e^{-t\lambda^2}d\lambda\\
&=&
 \frac{|L_{\infty}|C_{\rho}}{\pi^2}\cdot z\int_{-\infty}^{\infty}\frac{d\lambda}{\lambda^2+z^2}\\
 &=& \frac{|L_{\infty}|C_{\rho}}{\pi}.
\end{eqnarray*}
The same computation implies
\[
 L^{\prime}(U_1)(z)=\frac{2|L_{\infty}|C_{\rho}}{\pi}.
\]
Thus we have proved
\begin{prop}
Suppose $\rho|_{\Gamma_{\infty}}$ is nontrivial. Then we have
\[
 L^{\prime}(e^tU_0)(z-1)-L^{\prime}(U_1)(z)+L^{\prime}(e^tU_0)(z+1)=0.
\]
\end{prop}
In the following, we will compute
$2|L_{\infty}|R_{\rho}T^{\prime}(k_{j,t},0)$. Since these terms
appear only if $\rho|_{\Gamma_{\infty}}$ is trivial, we will assume
its triviality. Under the assumption, $R_{\rho}$ is explicitly computed as
\[
 R_{\rho}=\frac{\pi}{|L_{\infty}|},
\]
and we have
\[
 2|L_{\infty}|R_{\rho}T^{\prime}(k_{j,t},0)=2{\pi}T^{\prime}(k_{j,t},0).
\]
Let
\[
 SL_{2}({\mathbb C})\stackrel{\pi}\to G=PSL_{2}({\mathbb C})
\]
be the universal covering and $\tilde{K}$ (resp. $\tilde{M}$) the
inverse image of $K$ (resp. $M$), which is isomorphic to
$SU_{2}({\mathbb C})$ (resp. $U(1)$). $SL_{2}({\mathbb C})$ acts on the
polynomial ring ${\mathbb C}[X,Y]$ by the linear transformation, which is decomposed into irreducible representations:
\[
 {\mathbb C}[X,Y]=\oplus_{\mu}\xi_{\mu}.
\]
Here $\mu$ runs through nonnegative half integers and $\xi_{\mu}$ is the
space of homogeneous polynomial of degree $2\mu$. In particular
$\tau_{1}\otimes{\mathbb C}$ is isomorphic to $\xi_{1}$ as
$\tilde{K}$-modules. Let $D_{\tilde{K}}$ (resp.  $D_{\tilde{M}}$) be the set
of irreducible representations of $\tilde{K}$ (resp. $\tilde{M}$),
which will be parametrized by the set of nonnegative half integers
$\frac{1}{2}{\mathbb Z}_{\geq 0}$
(resp. half integers $\frac{1}{2}{\mathbb Z}$). Then the map
\[
 D_{K}\stackrel{\pi^{\dag}}\to D_{\tilde{K}}
\] 
and
\[
 D_{M}\stackrel{\pi^{\dag}}\to D_{\tilde{M}}
\]
may be identified with the natural inclusions:
\[
 {\mathbb Z}_{\geq 0} \hookrightarrow \frac{1}{2}{\mathbb Z}_{\geq 0},
\]
and
\[
 {\mathbb Z} \hookrightarrow \frac{1}{2}{\mathbb Z},
\]
respectively. \\

Now here is the Hoffmann's formula. (See the equation (8) and
(49) in \cite{Hoffmann}. See also \cite{Park}, $\S4$):
\begin{enumerate}
\item
\[
 T^{\prime}(k_{0,t},0)=\frac{1}{\pi}\{J^{(1)}_{0}(t)+\frac{1}{2\pi}{\rm p.v.}\int^{\infty}_{-\infty}J^{(2)}(0,i\lambda)(t)d\lambda\},
\]
\item 
\[
 T^{\prime}(k_{1,t},0)=\frac{1}{\pi}\sum_{\sigma=-1}^{1}\{J^{(1)}_{\sigma}(t)+\frac{1}{2\pi}{\rm p.v.}\int^{\infty}_{-\infty}J^{(2)}(\sigma,i\lambda)(t)d\lambda\},
\]
\end{enumerate}
where ${\rm p.v.}$ means Cauchy's principal value. We will explain
each term.
\begin{enumerate}
\item  $J^{(1)}_{\sigma}(t)$\\

Let $\psi$ be the di-gamma function:
\[
 \psi(z)=\frac{d}{dz}\log \Gamma (z)=\frac{\Gamma^{\prime}(z)}{\Gamma(z)}.
\]
It is known that it has a logarithmic growth as $|\lambda|\to\infty$:
\[
 |\psi(i\lambda)|\sim \log|\lambda|.
\]
We define a function
      $\Omega(\sigma,-i\lambda)$ for $\lambda\in{\mathbb R}$ as
\[
 \Omega(0,i\lambda)=2\psi(1)-(\psi(1+i\lambda)+\psi(1-i\lambda))
\]
and
\[
 \Omega(-1,i\lambda)=\Omega(1,i\lambda)\\
= 2\psi(1)-\frac{1}{2}(\psi(-i\lambda)+\psi(i\lambda)+\psi(2-i\lambda)+\psi(2+i\lambda)).
\]
{\bf Corollary} of p.96 in \cite{Hoffmann} and the computation of $\S4$
      of \cite{Park} shows
\[
 J^{(1)}_{\sigma}(t)=\frac{1}{2\pi}\int^{\infty}_{-\infty}\Omega(\sigma,-i\lambda)\hat{k}_{j,t}(\sigma, i\lambda)d\lambda.
\] 
Note that the integral is absolutely convergent because
      $\psi(z)+\psi(-z)$ is regular at $z=0$. We will compute
      $L^{\prime}(e^tJ^{(1)}_{0})(z)$ and
      $L^{\prime}(J^{(1)}_{\sigma})(z)$ for $\sigma=\pm 1$.
\begin{lm}
Let $z$ and $\alpha$ be positive numbers. Then we have the following
 identities.
\begin{enumerate}
\item
\[
 L^{\prime}(\int^{\infty}_{-\infty}\psi(\alpha+i\lambda)e^{-t\lambda^2}d\lambda)=L^{\prime}(\int^{\infty}_{-\infty}\psi(\alpha-i\lambda)e^{-t\lambda^2}d\lambda)=2\pi\psi(z+\alpha).
\]
\item
\[
 L^{\prime}(\int^{\infty}_{-\infty}\frac{e^{-t\lambda^2}}{\alpha+i\lambda}d\lambda)=L^{\prime}(\int^{\infty}_{-\infty}\frac{e^{-t\lambda^2}}{\alpha-i\lambda}d\lambda)=\frac{2\pi}{z+\alpha}.
\]
\item
\[
 L^{\prime}(\int^{\infty}_{-\infty}e^{-t\lambda^2}d\lambda)=2\pi.
\]
\end{enumerate}
\end{lm}
{\bf Proof.}
Changing the order of integration, we have
\begin{eqnarray*}
L^{\prime}(\int^{\infty}_{-\infty}\psi(\alpha+i\lambda)e^{-t\lambda^2}d\lambda)&=& 2z\int^{\infty}_{0}dt e^{-tz^2}\int^{\infty}_{-\infty}\psi(\alpha+i\lambda)e^{-t\lambda^2}d\lambda\\
&=& 2z \int^{\infty}_{-\infty}\frac{\psi(\alpha+i\lambda)}{\lambda^2+z^2}d\lambda.
\end{eqnarray*}
Note that $\psi(\alpha+i\lambda)$ is regular on ${\rm Im}\,\lambda<\alpha$ and  that $|\psi(\alpha+i\lambda)|\sim \log |\lambda|$ as $|\lambda|\to \infty$ (cf. \cite{Park}, p.22). Making a curvilinear integral along a contour which first goes along the real axis from $-\infty$ to $\infty$ then turns around the under semi-circle by the clock-wise direction, we will obtain
\begin{eqnarray*}
 \int^{\infty}_{-\infty}\frac{\psi(\alpha+i\lambda)}{\lambda^2+z^2}d\lambda&=&-2\pi i\cdot{\rm Res}_{\lambda=-iz}[\frac{\psi(\alpha+i\lambda)}{\lambda^2+z^2}d\lambda]\\
&=& \frac{\pi}{z}\psi(z+\alpha),
\end{eqnarray*}
which implies the first identity. The remaining identities can be proved by the
      same way. 
\begin{flushright}
$\Box$
\end{flushright}
Using {\bf Lemma 7.2}, we will obtain 
\[
 L^{\prime}(e^tJ^{(1)}_{0})(z)=2(\psi(1)-\psi(z+1)).
\]
Also, for a positive number $w$, {\bf Lemma 7.2} shows
\[
 L^{\prime}(\frac{1}{4\pi}\int^{\infty}_{-\infty}e^{-tz^2}\{\psi(w+i\lambda)+\psi(w-i\lambda)\}d\lambda)(z)=\psi(z+w).
\]
Making $w\to 0$, we obtain
\[
 L^{\prime}(\frac{1}{4\pi}\int^{\infty}_{-\infty}e^{-tz^2}\{\psi(i\lambda)+\psi(-i\lambda)\}d\lambda)(z)=\psi(z),
\]
which implies
\[
 L^{\prime}(J^{(1)}_{-1})(z)=L^{\prime}(J^{(1)}_{1})(z)=2\psi(1)-\psi(z)-\psi(z+2).
\]
Thus we have prove the following proposition.
\begin{prop}
\begin{enumerate}
\item
\[
  L^{\prime}(e^tJ^{(1)}_{0})(z)=2(\psi(1)-\psi(z+1)).
\]
\item
\[
 L^{\prime}(J^{(1)}_{-1})(z)=L^{\prime}(J^{(1)}_{1})(z)=2\psi(1)-\psi(z)-\psi(z+2).
\]
\end{enumerate}
In particular
\[
 L^{\prime}(e^tJ^{(1)}_{0})(z-1)+L^{\prime}(e^tJ^{(1)}_{0})(z+1)-\{L^{\prime}(J^{(1)}_{-1})(z)+L^{\prime}(J^{(1)}_{1})(z)\}=0.
\]
\end{prop}

\item $J^{(2)}(\sigma,i\lambda)(t)$\\

Using the logarithmic derivative of the Harish-Chandra's C-function $C(\sigma,\nu)$, Park
      computed the term as (see $\S4$ of \cite{Park}):
\[
 J^{(2)}(\sigma,i\lambda)(t)=-\hat{k}_{j,t}(\sigma, i\lambda)\cdot\frac{d}{d\nu}\log C (\sigma,\nu)|_{\nu=i\lambda}.
\]
\cite{Eguchi-Koizumi-Mamiuda} {\bf Theorem 8.2} and the functional
       equation:
\[
 \Gamma(z+1)=z\cdot\Gamma(z)
\]
      shows
\[
 C(\sigma,\nu)=
\left\{
\begin{array}{ccc}
\frac{1}{\nu} & if & \sigma=0\\
\nu & if & \sigma=\pm 1.
\end{array}
\right.
\]
Therefore $J^{(2)}(\sigma,i\lambda)(t)$ is computed as
\[
 J^{(2)}(\sigma,i\lambda)(t)=
\left\{
\begin{array}{ccc}
\frac{e^{-t(1+\lambda^2)}}{i\lambda} & if & \sigma=0\\
-\frac{e^{-t\lambda^2}}{i\lambda} & if & \sigma=\pm 1.
\end{array}
\right.
\]
Therefore 
\begin{eqnarray*}
\frac{1}{2\pi}{\rm p.v.}\int^{\infty}_{-\infty}e^tJ^{(2)}(0,i\lambda)(t)d\lambda&=& \frac{1}{2\pi}{\rm p.v.}\int^{\infty}_{-\infty}J^{(2)}(-1,i\lambda)(t)d\lambda\\
&=& \frac{1}{2\pi}{\rm p.v.}\int^{\infty}_{-\infty}J^{(2)}(1,i\lambda)(t)d\lambda=0.
\end{eqnarray*}

\end{enumerate}
Together with {\bf Corollary 7.1}, {\bf Proposition 7.2} and {\bf
Proposition 7.3}, these computation shows
\begin{prop}
Suppose $\rho|_{\Gamma_{\infty}}$ is trivial. Then
\[
 L^{\prime}(e^tU_{0})(z-1)-L^{\prime}(U_{1})(z)+L^{\prime}(e^tU_{0})(z+1)=0.
\]
\end{prop}


\section{Scattering terms}

In this section we assume that the restriction of $\rho$ to
$\Gamma_{\infty}$ is trivial. We first recall the basic facts of the
scattering matrix. (See \cite{EGM} and \cite{Sarnak-Wakayama})\\

Let $E_{\infty}(\varphi_{\sigma},\rho,s)$
be the constant term of the Fourier expansion of
$E(\varphi_{\sigma},\rho,s)$ at the cusp:
\[
 E_{\infty}(\varphi_{\sigma},\rho,s)(x)=\frac{1}{|L_{\infty}|}\int_{L_{\infty}}E(\varphi_{\sigma},\rho,s)(xn)dn.
\]
Then there is an intertwining operator called {\it the scattering matrix}:
\[
 {\mathcal H}(\sigma,\,s)\stackrel{C_{\rho,\sigma}(s)}\to {\mathcal H}(\sigma,\,s),
\]
which satisfies
\[
 E_{\infty}(\varphi_{\sigma},\rho,s)(x)=r(x)^{1+s}\varphi_{\sigma}(x)+r(x)^{1-s}(C_{\rho,\sigma}(s)\varphi_{\sigma})(x).
\]
In fact $C_{\rho,\sigma}(s)$ is a scalar if $\sigma$ is $0$ and is a $2\times
2$-matrix if $\sigma=\pm 1$. In the latter case, since
${\mathcal H}_{1,s}$ is not isomorphic to ${\mathcal H}_{-1,s}$, we know
(see \cite{Park} $\S3$):
\begin{equation}
 C_{\rho,\sigma}(s)=
\left(
\begin{array}{cc}
0 & C_{+}(\sigma, s)\\
C_{-}(\sigma, s) & 0
\end{array}
\right).
\end{equation}
The arugument of \cite{Friedman} $\S 3.7$ (see also \cite{EGM} $\S 6$
and \cite{Sarnak-Wakayama} $\S1$, $\S2$) shows the scattering matrix satisfies the following
properties:
\begin{enumerate}
\item
\[
 C_{\rho,\sigma}(0)=I.
\]
\item $C_{\rho,\sigma}(s)$ is absolutely convergent for ${\rm Re}\, s>1$ and is
       continued on the whole plane as a meromorphic function of order four.
\item It satisfies the functional equation:
\[
 C_{\rho,\sigma}(-s)\cdot C_{\rho,\sigma}(s)=I,
\] 
and its transpose conjugate $C_{\rho,\sigma}(s)^{*}$ is equal to $C_{\rho,\sigma}(\bar{s})$.
\end{enumerate}
In particular $C_{\rho,\sigma}(i\lambda)$ is a unitary matrix for
$\lambda\in {\mathbb R}$. \\

{\bf The scattering terms.}\\

Let us put
\[
 \Psi_{\rho,\sigma}(s)={\rm Tr}[C_{\rho,\sigma}(-s)\frac{d}{ds}C_{\rho,\sigma}(s)].
\]

Note that by definition we have
\begin{equation}
 \Psi_{\rho,\sigma}(s)=\Psi_{\rho,-\sigma}(s), \quad \sigma=-1,\,0,\,1.
\end{equation}

Using {\bf Fact 7.1} and the formula in p.16 of \cite{Sarnak-Wakayama}
together with (11),
we will see
\[
 {\mathcal S}_{0}(t)=\frac{1}{4\pi}\int^{\infty}_{-\infty}\Psi_{\rho,0}(i\lambda)e^{-t(1+\lambda^2)}d\lambda,
\]
and 
\[
 {\mathcal S}_{1}(t)=\frac{1}{4\pi}\int^{\infty}_{-\infty}\Psi_{\rho,0}(i\lambda)e^{-t(1+\lambda^2)}d\lambda+\frac{1}{2\pi}\int^{\infty}_{-\infty}\Psi_{\rho,1}(i\lambda)e^{-t\lambda^2}d\lambda.
\]
We put
\[
 S_{0}(t)={\mathcal S}_{0}(t),
\]
\[
 S_{1}(t)={\mathcal S}_{1}(t)-{\mathcal S}_{0}(t)=\frac{1}{2\pi}\int^{\infty}_{-\infty}\Psi_{\rho,1}(i\lambda)e^{-t\lambda^2}d\lambda.
\]
Following the argument of $\S6.4$ of \cite{EGM} we will obtain
\[
 \Psi_{\rho,0}(s)=c_{\rho,0}-\sum_{k}(\frac{1}{s-\alpha_k}-\frac{1}{s+\overline{\alpha_k}}),
\]
where $c_{\rho,0}$ is a constant and $\{\alpha_k\}_k$ is the set of
poles of $C_{\rho,0}(z)$. Let $\{\beta_l\}_l$ be the set of poles of
$\det C_{\rho,1}(s)$. Then the same argument as p.33 of \cite{Park} also
implies
\[
  \Psi_{\rho,1}(s)=c_{\rho,1}-\sum_{k}(\frac{1}{s-\beta_l}-\frac{1}{s+\overline{\beta_l}}),
\]
where $c_{\rho,1}$ is a constant. Since both $C_{\rho,0}(z)$ and $\det
C_{\rho,1}(s)$ are regular on the imaginary axis, real part of any
$\alpha_k$ and $\beta_l$ are nonzero. Thus
\[
 \Psi_{\rho,0}(i\lambda)d\lambda=\{c_{\rho,0}-\frac{1}{i}\sum_{k}(\frac{1}{\lambda+i\alpha_k}-\frac{1}{\lambda-i\overline{\alpha_k}})\}d\lambda
\]
and
\[
 \Psi_{\rho,1}(i\lambda)d\lambda=\{c_{\rho,1}-\frac{1}{i}\sum_{k}(\frac{1}{\lambda+i\beta_l}-\frac{1}{\lambda-i\overline{\beta_l}})\}d\lambda
\]
are 1-forms regular on the real axis.\\

For $z>0$ and a complex number $\alpha$ whose imaginary part
nonzero, we obtain by {\bf Lemma 7.1} 
\begin{eqnarray*}
 2z\int_{0}^{\infty}dte^{-tz^2}\int^{\infty}_{-\infty}\frac{e^{-t\lambda^2}}{\lambda-\alpha}d\lambda&=&2z\int^{\infty}_{-\infty}d\lambda \frac{1}{\lambda-\alpha}\int_{0}^{\infty}e^{-t(\lambda^2+z^2)}dt\\
&=& 2z\int^{\infty}_{-\infty}\frac{d\lambda}{(\lambda^2+z^2)(\lambda-\alpha)}\\
&=& \frac{2\pi i}{z-{\rm sgn}({\rm Im}\,\alpha)i\alpha},
\end{eqnarray*} 
which implies
\begin{eqnarray*}
 L^{\prime}(e^tS_{0})(z)&=&\frac{1}{2}\{c_{\rho,0}-\sum_{k}(\frac{1}{z+{\rm
 sgn}({\rm Re}\,\alpha_k)\alpha_k}-\frac{1}{z+{\rm
 sgn}({\rm Re}\,\alpha_k)\overline{\alpha_k}})\}\\
&=& \frac{1}{2}\{c_{\rho,0}-\sum_{{\rm Im}\,\alpha_k\neq 0}(\frac{1}{z+{\rm
 sgn}({\rm Re}\,\alpha_k)\alpha_k}-\frac{1}{z+{\rm
 sgn}({\rm Re}\,\alpha_k)\overline{\alpha_k}})\},
\end{eqnarray*}
and
\begin{eqnarray*}
 L^{\prime}(S_{1})(z)&=&c_{\rho,1}-\sum_{l}(\frac{1}{z+{\rm
 sgn}({\rm Re}\,\beta_l)\beta_l}-\frac{1}{z+{\rm
 sgn}({\rm Re}\,\beta_l)\overline{\beta_l}})\\
&=& c_{\rho,1}-\sum_{{\rm Im}\,\beta_l\neq 0}(\frac{1}{z+{\rm
 sgn}({\rm Re}\,\beta_l)\beta_l}-\frac{1}{z+{\rm
 sgn}({\rm Re}\,\beta_l)\overline{\beta_l}}).
\end{eqnarray*}
Thus we have proved the following proposition.
\begin{prop}
$L^{\prime}(e^tS_{0})(z)$ (resp. $L^{\prime}(S_{1})(z)$) is meromorphically continued on the
      whole plane with only simple poles whose residues are
      half integers (resp. integers). Moreover it is regular on the real axis.
\end{prop}

{\bf The threshold terms}\\

If $\sigma=\pm 1$, we know the trace of $C_{\rho,\sigma}(s)$
vanishes by (10) and the formula of p.16 of \cite{Sarnak-Wakayama} implies
\[
 {\mathcal T}_{0}(t)={\mathcal T}_{1}(t)=-\frac{1}{4}e^{-t}.
\]
Thus for $z>0$ we have
\begin{eqnarray*}
 L^{\prime}(e^t{\mathcal T}_{0})(z)&=&L^{\prime}(e^t{\mathcal T}_{1})(z)\\
&=& -\frac{1}{4}\cdot 2z\int^{\infty}_{0}e^{-tz^2}dt\\
&=& -\frac{1}{2z}.
\end{eqnarray*}
Now if we put
\[
 L^{\prime}_{0,sc}(z)=L^{\prime}(e^tS_{0})(z-1)+L^{\prime}(e^t{\mathcal T}_{0})(z-1),
\]
and
\[
 L^{\prime}_{1,sc}(z)=L^{\prime}(S_{1})(z),
\]
we have proved the following proposition.
\begin{prop}
$L^{\prime}_{0,sc}(z)$ and $L^{\prime}_{1,sc}(z)$ are continued on the
 whole plane as meromorphic functions with only simple poles whose
 residues are half integers and integers, respectively. Moreover we
 have
\[
 {\rm Res}_{z=0}L^{\prime}_{0,sc}(z)={\rm Res}_{z=2}L^{\prime}_{0,sc}(z)={\rm Res}_{z=0}L^{\prime}_{1,sc}(z)=0.
\]
\end{prop}

\section{A proof the theorem }

Now we are ready to prove {\bf Theorem 1.1}. Up to the previous section we
have seen the logarithmic derivative of the Ruelle L-function
$R_{\rho}(z)$ is
continued on the whole plane as a meromorphic function which has only
simple poles whose residues are at most half integers. Thus
$R_{\rho}(z)^2$ is meromorphically continued on the whole plane. But if
${\rho}|_{\Gamma_\infty}$ is nontrivial, since the contribution of the
scattering terms does not exist, we find all the residues are
integers. Thus $R_{\rho}(z)$ is
meromorphically continued on ${\mathbb C}$ itself.\\

By the Selberg trace formula and (7) we will find the set of poles of $\frac{d}{dz}\log R_{\rho}(z)$ is a union of two subsets, ${\rm P}_{reg}$
and ${\rm P}_{sing}$. Here ${\rm P}_{reg}$ (resp. ${\rm P}_{sing}$) is the set of eigenvalues of
$\Delta$ on $L^2(X,\,\Omega^{\cdot}(\rho))$ (resp. the poles of the
derivative of the Laplace transforms of the unipotent orbital integrals
an the scattering terms). Note that the derivative of the Laplace
transform of the identity orbital integrals does
not yield any pole. {\bf Proposition 7.4} and {\bf Proposition 8.2}
implies $0$ is not contained in ${\rm P}_{sing}$. Now we obtain
\begin{equation}
 {\rm ord}_{z=0}R_{\rho}(z)=2(2\beta_0(\rho)_{(2)}-\beta_1(\rho)_{(2)}).
\end{equation}
by {\bf Lemma 3.1}, {\bf Lemma 3.2} and the Selberg trace formula.
Together with {\bf Lemma 2.1} and {\bf Lemma 2.2}, (12) shows
{\bf Theorem 1.1}.

\begin{flushright}
$\Box$
\end{flushright}

Comparing the Riemann's zeta function, we may consider ${\rm
P}_{sing}$ or ${\rm P}_{reg}$
corresponds to the set of trivial or of 
essential zeros, respectively. As we have seen in $\S3$, except for
finitely many elements, it is
contained in
\[
 \{z\in{\mathbb C}\,|\, {\rm Re}\,z \in \{-1,\,0,\,1\}\}.
\]
Thus we may say the Ruelle L-function satisfies {\it the Riemann
hypothesis.} Note that if $\rho|_{\Gamma_{\infty}}$ is nontrivial,
{\bf Proposition 7.2} shows ${\rm P}_{sing}$ is empty.


\vspace{10mm}
\begin{flushright}
Address : Department of Mathematics and Informatics\\
Faculty of Science\\
Chiba University\\
1-33 Yayoi-cho Inage-ku\\
Chiba 263-8522, Japan \\
e-mail address : sugiyama@math.s.chiba-u.ac.jp
\end{flushright}
\end{document}